\documentclass{cimart}

\usepackage{enumitem}
\usepackage{mathrsfs}

\newcommand{\R}{\mathbb{R}}
\newcommand{\N}{\mathbb{N}}

\newtheorem*{examples}{Examples}

\title{Weighted norm inequalities for integral transforms with splitting kernels}

\author{Alberto Debernardi Pinos}

\authorinfo{Departament de Matem\`atiques, Universitat Aut\`onoma de Barcelona, Spain}{adebernardipinos@gmail.com}

\abstract{%
	We obtain necessary and sufficient conditions on weights for a wide class of integral transforms to be bounded between weighted $L^p-L^q$ spaces, with $1 \le p \le q \le \infty$. The kernels $K(x,y)$ of such transforms are only assumed to satisfy upper bounds given by products of two functions, one in each variable.
	
	The obtained results are applicable to several transforms, some of which are included here as particular examples. Some of the new results derived here are the characterization of weights for the boundedness of the $\mathscr{H}_\alpha$ (or Struve) transform in the case $\alpha>\frac{1}{2}$, or the characterization of power weights for which the Laplace transform is bounded in the limiting cases $p=1$ or $q=\infty$.
}



\VOLUME{32}
\ISSUE{3}
\firstpage{241}
\DOI{https://doi.org/10.46298/cm.12756}

\begin{document}
	
\section{Introduction}

Given a linear operator $T$, studying its boundedness between different pairs of function spaces is an important problem in analysis. If furthermore $T$ is an integral operator, often necessary and sufficient conditions for the boundedness of such an operator between two spaces may be expressed in terms of the finiteness of appropriate integral expressions.

In this paper, we study integral transforms of the type
\begin{equation}
	\label{trans}
	Tf(y)=\int_0^\infty f(x)K(x,y)\, dx, \qquad y\in \R_+:=(0,\infty),
\end{equation}
with the kernel $K$ satisfying an upper estimate of the form
\begin{equation}
	\label{kernelest}
	|K(x,y)|\leq  \begin{cases}
		C_1s_1(x)w_1(y), & x\in (0,\varphi(y)),\\
		C_2s_2(x)w_2(y), &x\in (\varphi(y),\infty),
	\end{cases}
\end{equation}
where $C_j>0$, $\varphi:(0,\infty)\to (0,\infty)$ is a bijective function of class $C^1$, and $s_1, s_2, w_1, w_2$ are nonnegative functions. In this case we call $T$ an \textit{integral transform with splitting kernel}, motivated by the fact that the upper estimate \eqref{kernelest} splits into a product of a function of $x$ and a function of $y$ (we emphasize that the kernel need not split this way, but only the required upper estimate).

Our goal is to study the boundedness of integral transforms with splitting kernels between pairs of weighted Lebesgue spaces. More specifically, we aim to obtain necessary and sufficient conditions on the weights (positive locally integrable functions) $u$ and $v$ for the weighted norm inequality
\begin{equation}
	\label{wni}
	\Vert Tf\Vert_{L_q^u}\leq C_{T,p,q,u,v} \Vert f\Vert_{L_p^v}
\end{equation}
to hold with $1\leq p\leq q\leq \infty$, where $\|\cdot \|_{L_p^v}$ denotes the weighted Lebesgue norm given by
\[
\Vert f\Vert_{L_p^v}=\bigg(\int_0^\infty v(x)^p|f(x)|^p\, dx\bigg)^\frac{1}{p},
\]
with the obvious modification for $p=\infty$. We emphasize that we also consider the cases $p=1$ and $q=\infty$, which are often excluded in such problems. If $v\equiv 1$, we simply write $\| f\|_{L_p^v}=\|f\|_{L_p}$. 
In some cases, we have to consider the $L_p$ (or $L_p^v$) norm of a function over a measurable subset $E\subset \R_+$, for which we denote
\[
\|f\|_{L_p^v(E)}:=\bigg( \int_E v(x)^p |f(x)|^p\, dx\bigg)^\frac{1}{p},
\]
and similarly for $\|f\|_{L_p(E)}$. Finally, a function $f$ is said to be locally integrable ($f\in L_1^{\textrm{loc}}$) if $\|f\|_{L_1(E)}<\infty$ for every compact set $E\subset \R_+$. 

In the sequel we will assume that the function $f$ in the transform \eqref{trans} is such that $f(x)K(x,y)\in L_1^{\textrm{loc}}$ for every $y>0$, and that the corresponding integral transform exists pointwise as an improper integral, i.e., for all $y>0$, the limit
\[
\lim_{N\to \infty}\int_0^N f(x)K(x,y)\, dx
\]
exists and is finite.

Different types of integral transforms require different techniques for obtaining sharp weighted norm inequalities between Lebesgue spaces. We begin by mentioning weighted Hardy-type inequalities, see, e.g., \cite{hardyineq,MuckenhouptHardy,talenti,tomaselli}. These inequalities usually serve as a basis for deriving corresponding weighted norm inequalities for other integral operators (in particular, the Hardy-type inequalities in Subsection~\ref{SUBSEChardy} below are the main tool we use in this paper). 

Another related example is the Riemann-Liouville fractional integral, which is a direct generalization of the Hardy operator. In \cite{MRS}, weighted norm inequalities are characterized in terms of the weights (which resemble those conditions in Hardy's inequalities) for the Riemann-Liouville integral of order larger than one. For slightly more general operators that include the Riemann-Liouville integral of any order, the full characterization of such inequalities was obtained slightly later in \cite{BK}.

We also refer to the recent paper \cite{liflyand}, where necessary and sufficient conditions for weighted inequalities (with power weights) for Hausdorff operators are studied. The prototype of Hausdorff operator is precisely the Hardy operator. However, Hausdorff operators are rather general and include several well-known examples, as for instance, the aforementioned Riemann-Liouville integrals \cite{BK,MRS}.

In relation to integral transforms of Fourier type, sharp inequalities can be obtained through an approach based on rearrangement inequalities, which, in essence, have an underlying interpolation argument. Such inequalities are not intrinsic to the particular characteristics of the Fourier transform, but to the more general operators of type $(2,2)$ and $(1,\infty)$ (i.e., so that they are bounded in $L_2$, and also from $L_1$ to $L_\infty$). Fourier-type transforms, as for instance, the sine, cosine, or Hankel transforms, and many others such as the Laplace transform satisfy this kind of estimates.

Necessary and sufficient conditions for such inequalities were first obtained in \cite{heinig} for the cosine transform, and ultimately for the Fourier transform. We also mention the paper  \cite{RonCan}, where sufficient conditions for inequalities with power weights are derived for certain Hankel-type transforms, the paper \cite{decarli}, where inequalities with power weights were characterized for the Hankel transform (i.e., the Fourier transform of radial functions), and the recent work \cite{GLTPitt}, where a unified approach for Fourier-type transforms was developed, and sharp conditions for sine transform inequalities involving power weights were obtained for the first time (surprisingly enough, given the considerable amount of literature in the topic). In the last reference, necessary and sufficient conditions involving general weights were also derived.

The case of the Laplace transform, although being of types $(2,2)$ and $(1,\infty)$ as the Fourier transform, needs different treatment if one desires to obtain sharp conditions for norm inequalities. This is because, as we discuss in more detail in Subsection~\ref{SUBSEClaplace}, the aforementioned interpolation type estimates are based on decreasing rearrangements of the weights. However, given the fast decay of the kernel $K(x,y)=e^{-xy}$ of the Laplace transform, direct estimates are sometimes unavoidable, especially in the case where the involved weights have a decreasing rearrangement that is identically infinity. This is illustrated in the two main results of \cite{BloomLaplace}, which derives norm inequalities for the Laplace transform using these two approaches (see also \cite{SU}; it is worth noting that direct estimates are based on the self-adjointness of the operator defined by the Laplace transform).

In the sequel we use an approach for deriving inequalities of the form \eqref{wni} based on direct estimates combined with Hardy-type inequalities. More precisely, we will see in Section~\ref{SECwni} that an estimate of the type \eqref{kernelest} allows a rather direct application of Hardy's inequalities, which in turn yields sufficient conditions for \eqref{wni} to hold in terms of the weights and of the functions $s_j,w_j$, $j=1,2$, that appear in \eqref{kernelest}.

On the other hand, for deriving necessary conditions for \eqref{wni} to hold we require that the kernel $K(x,y)$ is positive, since necessary conditions are usually based on lower estimates for $K(x,y)$ similar to \eqref{kernelest}. In particular, if a reverse estimate of the form
\begin{equation*}
	K(x,y)\geq  \begin{cases}
		c_1s_1(x)w_1(y), & x\in (0,\varphi(y)),\\
		c_2s_2(x)w_2(y), &x\in (\varphi(y),\infty),
	\end{cases}
\end{equation*}
holds for some $c_1,c_2>0$, the necessary and sufficient conditions for \eqref{wni} coincide (see Corollary~\ref{CORiff}).

Applying the obtained results to the $\mathscr{H}_\alpha$ transform (cf. Subsection~\ref{SUBSEChtrans} below), we are able to characterize the inequality $\| \mathscr{H}_{\alpha} f\|_{L_q^u}\leq C_{p,q,\alpha,u,v}\| f\|_{L_p^v}$ by the condition
\[
\sup_{t>0}  \bigg( \int_0^\infty \bigg(  \frac{x^{\alpha+\frac{3}{2}}}{t^{-2}+x^2}\bigg)^{q} u(x)^{q}\, dx\bigg)\bigg( \int_0^\infty \bigg(  \frac{x^{\alpha+\frac{3}{2}}}{t^{2}+x^2}\bigg)^{p'} v(x)^{-p'}\, dx\bigg)<\infty,
\]
for $1<p\leq q<\infty$ and $\alpha>\frac{1}{2}$ in Theorem~\ref{CORgluningstruve}. We also characterize the corresponding inequality for the cases $p=1$ or $q=\infty$ in Theorem~\ref{THMiffstruve}.

We can also characterize inequalities involving the Laplace transform (which we denote by $\mathcal{L}$) in the particular case where the weights are power functions. More precisely, we show that for $1\leq p\leq q\leq \infty$, with $(p,q)\neq (1,\infty)$, the inequality  $\Vert x^{-\beta}\mathcal{L}f\Vert_{L_q}\lesssim \Vert x^\gamma f\Vert_{L_p}$
holds if and only if
\[
\beta<\frac{1}{q},\quad \text{and}\quad \beta =\gamma+\frac{1}{q}-\frac{1}{p'},
\]
see Corollary~\ref{CORlaplacepower}. Furthermore, we also show that the inequality $\| x^{-\beta}\mathcal{L}f\|_{L_\infty}\lesssim \| x^\gamma f\|_{L_{1}}$ holds if and only if $\beta=\gamma=0$. Such necessary and sufficient conditions were obtained in \cite{RonCan,SU} in the case $1<p\leq q<\infty$. It is worth emphasizing that, although the kernel $K(x,y)=e^{-xy}$ of the Laplace transform does not satisfy a sharp estimate of the form \eqref{kernelest}, it is still possible to characterize the inequalities involving power functions (in fact, our results do not allow for such a characterization with general weights).

The outline of the paper is as follows. In Section~\ref{SECexamples0} we list several operators that fall into the class of integral transforms with splitting kernel, providing explicit estimates of the form \eqref{kernelest}. In Section~\ref{SECprelim} we discuss some properties of the $\mathscr{H}_\alpha$ transform, and we also obtain Hardy-type inequalities adapted to integral transforms with splitting kernel. Section~\ref{SECwni} is devoted to obtaining necessary and sufficient conditions for inequalities of the form \eqref{trans} to hold. We also give corresponding ``gluing lemmas'', which allow to write Hardy-type conditions jointly in some particular cases. Finally, all these results are applied to specific transforms in Section~\ref{SECexamples}, along with further discussion.

\medskip

\section{Examples of integral transforms with splitting kernels}\label{SECexamples0}

Let us give examples of integral transforms of the form \eqref{trans}, whose kernel satisfies an estimate of the type \eqref{kernelest}. For each of these transforms, we give explicit expressions for the functions $s_j, w_j$, $j=1,2$, and $\varphi$.

It is important to emphasize that we are not requiring the estimate \eqref{kernelest} to be sharp. This means that several classical integral operators fall under the scope of our work, although some estimates may be rather rough. As we will see in Section~\ref{SECwni}, the sharpness of the obtained results will be directly related to the sharpness of the estimate \eqref{kernelest}, allowing characterizations of the weights $u$ and $v$ in \eqref{wni} whenever the estimate \eqref{kernelest} is actually an equivalence.

Before proceeding further, we introduce some notation. By $A\lesssim B$ (resp. $A\gtrsim B$) we mean that there exists a constant $C>0$ not depending on essential quantities such that $A\leq CB$ (resp. $A\geq CB$). If $A\lesssim B$ and $A\gtrsim B$ simultaneously, we write $A\asymp B$.

Some examples of integral transforms with splitting kernels are the following.

\begin{examples}
	\begin{itemize}[leftmargin=*]
		\item The Hardy operator 
		\[
		Hf(y)=\int_0^y f(x)\, dx
		\]
		is of the form \eqref{trans} with $\varphi(y)=y$,  $s_1=w_1\equiv 1$ and $s_2=w_2\equiv 0$.
		\item The Bellman operator
		\[
		Bf(y)=\int_y^\infty f(x)\frac{dx}{x}
		\]
		is of the form \eqref{trans} with $\varphi(y)=y$, $s_1=w_1\equiv 0$, $w_2\equiv 1$ and $s_2(x)=1/x$.
		\item The Riemann-Liouville integral
		\[
		I_\alpha f(y)=\int_0^y f(x)(y-x)^{\alpha-1}\, dx,\qquad \alpha\in (0,1).
		\]
		In this case,  $K(x,y)=(y-x)^{\alpha-1}\chi_{(0,y)}(x)$, and
		\[
		(y-x)^{\alpha-1}\chi_{(0,y)}(x)\leq \begin{cases}
			y^{\alpha-1},&\text{if }x<y,\\
			0,&\text{if }x>y,
		\end{cases}
		\]
		i.e., $\varphi(y)=y$, $s_1(t)=1$, $w_1(t)=t^{\alpha-1}$, and $s_2=w_2\equiv 0$. Of course, this estimate is far from being sharp. Nevertheless, this illustrates the mildness of the requirement of satisfying an estimate of the form \eqref{kernelest}.
		\item The sine transform 
		\[
		\widehat{f}_{\sin}(y)=\int_0^\infty f(x)\sin( xy)\, dx.
		\]
		Since $|\sin xy|\leq \min\{xy, 1\}$ for $x,y>0$, in this case $\varphi(y)=1/y$, $s_1(t)=w_1(t)=t$, and $s_2=w_2\equiv 1$. Note, however, that the estimate $|\sin xy|\leq 1$ is rather rough. This will translate into results that are not sharp.
		
		\item The $\mathscr{H}_\alpha$ (or Struve) transform \cite{RonCan,titchmarsh}
		\begin{equation*}
			\mathscr{H}_\alpha(y)=\int_0^\infty (xy)^{\frac{1}{2}} f(x)\mathbf{H}_\alpha(xy)\, dx, \qquad \alpha>-\frac{1}{2} \label{defHtrans},
		\end{equation*}
		where $\mathbf{H}_\alpha$ is the Struve function of order $\alpha$, cf. \cite{EMOT,watBesselfn}.  	The definition and further properties of $\mathscr{H}_\alpha$ and $\mathbf{H}_\alpha$ are given in Subsection~\ref{SUBSEChtrans} below. We will see in \eqref{Hest} that the kernel of the $\mathscr{H}_\alpha$ transform satisfies
		\begin{equation}
			\label{EQeststruvesmallalpha}
			(xy)^{\frac{1}{2}} |\mathbf{H}_\alpha(xy)|\lesssim \begin{cases}
				(xy)^{\frac{3}{2}+\alpha},&\text{if } x<\frac{1}{y},\\
				1,&\text{if }x>\frac{1}{y},
			\end{cases} \qquad \text{for }-\frac{1}{2}<\alpha<\frac{1}{2},
		\end{equation}
		and
		\begin{equation}
			\label{EQeststruvebigalpha}
			(xy)^{\frac{1}{2}} \mathbf{H}_\alpha(xy)\lesssim \begin{cases}
				(xy)^{\frac{3}{2}+\alpha},&\text{if } x<\frac{1}{y},\\
				(xy)^{-\frac{1}{2}+\alpha},&\text{if }x>\frac{1}{y},
			\end{cases} \qquad \text{for }\alpha\geq\frac{1}{2},
		\end{equation}
		(for $\alpha\geq \frac{1}{2}$, $\mathbf{H}_\alpha$ is nonnegative). Hence, $\mathbf{H}_\alpha$ satisfies the estimate \eqref{kernelest} with $\varphi(y)=\frac{1}{y}$, $s_1(t)=w_1(t)=t^{\frac{3}{2}+\alpha}$ and $s_2=w_2\equiv 1$ for $-\frac{1}{2}<\alpha<\frac{1}{2}$, or $s_2(t)=w_2(t)=t^{-\frac{1}{2}+\alpha}$ for $\alpha\geq \frac{1}{2}$.
		%
		\item The (generalized) Stieltjes transform 
		\[
		S_\lambda f(y)=\int_0^\infty \frac{f(x)}{(x+y)^\lambda}\, dx, \qquad \lambda >0.
		\]
		Since
		\begin{equation}
			\label{EQstieltjesestimate}
			\frac{1}{(x+y)^{\lambda}}\asymp \begin{cases}
				y^{-\lambda}, &\text{if }x< y,\\
				x^{-\lambda},&\text{if }x>y,
			\end{cases}
		\end{equation}
		\eqref{kernelest} holds with $\varphi(y)=y$, $s_1(t)=w_2(t)=1$, and $w_1(t)=s_2(t)=t^{-\lambda}$. As we will note, the fact that such an estimate is sharp (it is indeed an equivalence) translates into sharp conditions on the weights that guarantee corresponding norm inequalities (cf. Subsection~\ref{SUBSECstieltjes} below). 
		
		\item The Laplace transform 
		\[
		\mathcal{L}f(y)=\int_0^\infty f(x)e^{-xy}\, dx.
		\]
		Note that for any $n\in \N$, one has the estimate
		\begin{equation}
			\label{EQestexponential}
			e^{-xy}\leq \begin{cases}
				1,&\text{if }x<\frac{1}{y},\\
				n!(xy)^{-n},&\text{if }x>\frac{1}{y}.
			\end{cases}
		\end{equation}
		Thus, for the Laplace transform, the estimate \eqref{kernelest} holds with $\varphi(y)=\frac{1}{y}$, $s_1=w_1\equiv 1$, and $s_2(t)=w_2(t)=t^{-n}$ for any choice of $n$. 
	\end{itemize}
\end{examples}

There are many other examples of integral transforms of the form \eqref{trans} whose kernel satisfies \eqref{kernelest}, as for instance the Meijer (or Macdonald) transform \cite[Ch. 23]{Z} (with $K(x,y)=(xy)^\frac{1}{2}K_\nu(xy)$, where $K_\nu$ is the Bessel function of the second kind), or the general integral transforms whose kernel a Fox $H$-function (called $H$-transforms; see \cite{KS} for a comprehensive description of the theory). In fact, the authors from \cite{BJ} use a similar approach as the one we consider here in order to derive weighted norm inequalities for $H$-transforms. 

An example of an integral transform of the type \eqref{trans} whose kernel does not satisfy an estimate of the form \eqref{kernelest} is the Mellin transform, with $K(x,y)=x^{y-1}$.

\medskip

\section{Preliminaries}\label{SECprelim}

\subsection{The Struve function and the \texorpdfstring{$\mathscr{H}_\alpha$} \text{} transform}\label{SUBSEChtrans}
The Struve function $\mathbf{H}_\alpha$ is defined by the series
\begin{equation}\label{Hseries}
	\mathbf{H}_\alpha(x)=\bigg(\frac{x}{2}\bigg)^{\alpha+1}\sum_{k=0}^\infty \frac{(-1)^k(x/2)^{2k}}{\Gamma(k+3/2)\Gamma(k+\alpha+3/2)},
\end{equation}
see, e.g., \cite[\S 7.5.4]{EMOT} or \cite[\S 10.4]{watBesselfn}. It is a continuous function and it is related to the Bessel function $J_\alpha$ in the following way: $\mathbf{H}_\alpha$ is the solution of the non-homogeneous Bessel differential equation
\begin{equation}
	\label{EDO}
	x^2\frac{d^2f}{dx^2}+x\frac{df}{dx}+(x^2-\alpha^2)f=\frac{4(x/2)^{\alpha+1}}{\sqrt{\pi}\Gamma(\alpha+1/2)},
\end{equation}
with the terms corresponding to the solution of the homogeneous equation equal to zero. On the other hand, $J_\alpha$ is the solution of the homogeneous differential equation corresponding to \eqref{EDO} that is bounded at the origin for nonnegative $\alpha$. 

For $\alpha>-\frac{1}{2}$, the equivalence
\begin{equation*}
	\label{EQequivxsmallstruve}
	\mathbf{H}_\alpha(x) \asymp x^{\alpha+1}, \qquad x\leq 1,
\end{equation*}
holds. Indeed, in view of \eqref{Hseries}, we only need to show that for $x\leq 1$,
\[
\sum_{k=0}^\infty \frac{(-1)^k(x/2)^{2k}}{\Gamma(k+3/2)\Gamma(k+\alpha+3/2)}\asymp 1.
\]
On the one hand, the latter series is absolutely convergent for $x\leq 1$, and thus bounded from above. On the other hand,
\[
\sum_{k=0}^\infty \frac{(-1)^k(x/2)^{2k}}{\Gamma(k+3/2)\Gamma(k+\alpha+3/2)}\geq \frac{1}{\Gamma(3/2)\Gamma(\alpha+3/2)}\bigg(1-\frac{x^2}{10(\alpha+5/2)}\bigg) \asymp 1, \qquad x\leq 1.
\]

For large $x$, we have the following asymptotic expansion \cite[p. 332]{watBesselfn},
\begin{align}
	\mathbf{H}_\alpha(x)&= 	\bigg(\frac{\pi x}{2}\bigg)^{-\frac{1}{2}}\sin\Big(x-\frac{\alpha \pi}{2}-\frac{\pi}{4}\Big)+\frac{(x/2)^{\alpha-1}}{\Gamma(\alpha+1/2)\Gamma(1/2)}(1+O(x^{-2})), \label{eqEstimateStruveLarge}
\end{align}
from which we can deduce
\begin{align*}
	|\mathbf{H}_\alpha(x)|&\lesssim x^{-\frac{1}{2}}+x^{\alpha-1} \asymp x^{\max\{-\frac{1}{2},\, \alpha-1\}} =\begin{cases} x^{-\frac{1}{2}}, & \text{if }\alpha<\frac{1}{2}, \\
		x^{\alpha-1}, &\text{if }\alpha\geq \frac{1}{2},\end{cases}\qquad x\geq 1.
\end{align*}
Combining these estimates, we obtain 
\begin{equation}
	\label{Hest}
	|\mathbf{H}_\alpha(x)|
	\lesssim \begin{cases}\min\{x^{\alpha+1},x^{-\frac{1}{2}}\}, &\alpha<\frac{1}{2},\\
		\min\{x^{\alpha+1},x^{\alpha-1}\}, &\alpha\geq \frac{1}{2},
	\end{cases}
\end{equation}
which obviously implies \eqref{EQeststruvesmallalpha} and \eqref{EQeststruvebigalpha}.

\begin{remark}
	\label{remarkHbigalpha}
	For $\alpha\geq \frac{1}{2}$ and $x>0$,  $\mathbf{H}_\alpha(x)$ is nonnegative \cite[p. 337]{watBesselfn}. Moreover, it follows from  \eqref{eqEstimateStruveLarge} that if $\alpha>\frac{1}{2}$, then there is $x_0>1$ such that
	\begin{equation*}
		\label{eqEstimateHspecial}
		\mathbf{H}_\alpha(x) \asymp x^{\alpha-1}, \qquad x>x_0.
	\end{equation*}
	Thus,
	\begin{equation*}
		\label{EQequivstruve}
		\mathbf{H}_\alpha(x)\asymp \min\{x^{\alpha+1},x^{\alpha-1}\},\qquad \text{for } \alpha>\frac{1}{2}.
	\end{equation*}
\end{remark}

\subsection{Hardy's inequalities}\label{SUBSEChardy}
The main tools we use in order to obtain sufficient conditions for \eqref{wni} are Hardy's inequalities \cite{hardyineq,MuckenhouptHardy}. Let us introduce the following notation before. For a measurable function $g$ defined on $\R_+$ and a measurable set $E\subset \R_+$, we put
\[
I_E(g)= \int_E |g(x)|\,dx.
\]
Also, for $1\leq p\leq \infty$, we define $p'$ to be the H\"older conjugate of $p$, i.e., so that $\frac{1}{p}+\frac{1}{p'}=1$ holds.

\begin{lemma}\textnormal{\cite{hardyineq}}\label{LEMhardylema}
	Let $1\leq p \leq q \leq \infty$. For a pair of weights $u$ and $v$, there exists $B>0$ such that the inequality
	\[
	\|  I_{(0,\cdot)}(g) \|_{L_q^u}\leq B \| g\|_{L_p^v} 
	\]
	holds for every measurable $g$ if and only if
	\begin{equation*}
		\label{EQhardycond01}
		\sup_{r>0}\| u\|_{L_q(r,\infty)} \| v^{-1}\|_{L_{p'}(0,r)} <\infty.
	\end{equation*}
	Also, there exists $B>0$ such that the inequality
	\[
	\|  I_{(\cdot,\infty)}(g) \|_{L_q^u}\leq B \| g\|_{L_p^v} 
	\]
	holds for every measurable $g$ if and only if
	\begin{equation*}
		\label{EQhardycond02}
		\sup_{r>0}\| u\|_{L_q(0,r)} \| v^{-1}\|_{L_{p'}(r,\infty)} <\infty.
	\end{equation*}
\end{lemma}		

We will need generalized forms of Hardy's inequalities that involve changes of variables. More precisely, instead of considering integrals of the form
\[
I_{(0,y)}(g), \qquad I_{(y,\infty)}(g),
\] 
we consider integrals of the form
\[
I_{(0,\varphi(y))}(g),\qquad I_{(\varphi(y),\infty)}(g).
\]
For a set $E\subset \R_+$ and a function $\psi: \R_+\to \R_+$, denote by $\psi(E)$ the image of $E$ under~$\psi$.

\begin{lemma}
	\label{LEMhardy-ineq-changevariables}
	Let $1\leq p\leq q \leq\infty$, and let  $\varphi:\R_+\to \R_+$ be a $C^1$ bijective function. For a pair of weights $u$ and $v$, there exists $B>0$  such that the inequality
	\begin{equation}
		\| I_{(0,\varphi(\cdot))}(g) \|_{L_q^u}\leq B \| g\|_{L_p^v} 
		\label{EQhardygeneral1}
	\end{equation}
	holds for any measurable $g$ if and only if 
	\begin{equation*}
		\label{EQhardycond1}
		\sup_{r>0}\| u\|_{L_q(\varphi^{-1}(r,\infty))} \| v^{-1}\|_{L_{p'}(0,r)}<\infty.
	\end{equation*}
	Also, there exists $B>0$ such that the inequality
	\begin{equation}
		\label{EQhardygeneral2}
		\|  I_{(\varphi(\cdot),\infty)}(g) \|_{L_q^u}\leq B \| g\|_{L_p^v} 
	\end{equation}
	holds for any measurable $g$ if and only if 
	\begin{equation*}
		\label{EQhardycond2}
		\sup_{r>0}\| u\|_{L_q(\varphi^{-1}(0,r))} \| v^{-1}\|_{L_{p'}(0,r)}<\infty.
	\end{equation*}
\end{lemma}
\begin{proof}
	We only prove the first part (the one corresponding to \eqref{EQhardygeneral1}), since the second one is analogous. We distinguish between the cases $q<\infty$ and $q=\infty$. For the case $q<\infty$, applying the change of variables $y= \varphi^{-1}(s)$ on the left hand side of \eqref{EQhardygeneral1}, we get
	\[
	\|  I_{(0,\varphi(\cdot))}(g) \|_{L_q^u}=\bigg(\int_0^\infty  \frac{u(\varphi^{-1}(s))^q}{\varphi'(\varphi^{-1}(s))} I_{(0,s)}(g)^q ds\bigg)^\frac{1}{q},
	\]
	and by the classical Hardy's inequality (Lemma~\ref{LEMhardylema}), we have
	\[
	\bigg(\int_0^\infty \frac{u(\varphi^{-1}(s))^q}{\varphi'(\varphi^{-1}(s))}I_{(0,s)}(g)^q  ds\bigg)^\frac{1}{q}\leq B \| g\|_{L_p^v} 
	\]
	if and only if
	\[
	\sup_{r>0}\bigg(\int_r^\infty \frac{u(\varphi^{-1}(s))^q}{\varphi'(\varphi^{-1}(s))}\, ds \bigg)^\frac{1}{q} \| v^{-1}\|_{L_{p'}(0,r)}<\infty.
	\]
	The substitution $x=\varphi^{-1}(s)$ shows that the last condition is equivalent to
	\[
	\sup_{r>0} \bigg(\int_{x\in \varphi^{-1}(r,\infty)} u(x)^q\, dx \bigg)^\frac{1}{q} \| v^{-1}\|_{L_{p'}(0,r)}=\sup_{r>0}\| u\|_{L_q(\varphi^{-1}(r,\infty))} \| v^{-1}\|_{L_{p'}(0,r)}<\infty.
	\]
	Let us now prove the case $q=\infty$. It is clear that
	\begin{align*}
		\|  I_{(0,\varphi(\cdot))}(g) \|_{L_\infty^u}&=\sup_{y\in \R_+}  u(y)\int_0^{\varphi(y)} |g(x)|\, dx=\sup_{y\in \R_+}  u(\varphi^{-1}(y))\int_{0}^y |g(x)|\, dx\\
		&= \|u(\varphi^{-1} (\cdot))I_{(0,\cdot)}(g) \|_{L_\infty}.
	\end{align*}
	By the classical Hardy's inequality, we have that
	\[
	\| u(\varphi^{-1} (\cdot))I_{(0,\cdot)}(g) \|_{L_\infty}\leq B \| g\|_{L_p^v} 
	\]
	if and only if
	\[
	\sup_{r>0}\|  u(\varphi^{-1} (\cdot))\|_{L_\infty(r,\infty)} \| v^{-1}\|_{L_{p'}(0,r)} <\infty,
	\]
	or equivalently, if and only if
	\[
	\sup_{r>0}\| u\|_{L_\infty (\varphi^{-1}(r,\infty))}\| v^{-1}\|_{L_{p'}(0,r)} <\infty.\qedhere
	\]
\end{proof}
Note that the classical Hardy's inequalities (Lemma~\ref{LEMhardylema}) correspond to Lemma~\ref{LEMhardy-ineq-changevariables} with $\varphi(r)=r$.

\medskip

\section{Weighted norm inequalities: necessary and sufficient conditions}\label{SECwni}
%
In this section we give necessary and sufficient conditions that guarantee that the inequality $\Vert Tf\Vert_{L_q^u}\lesssim \Vert f\Vert_{L_p^v}$ holds. We will also give two ``gluing lemmas'', aimed at writing jointly Hardy-type conditions.

\subsection{Sufficient conditions}

A  rather straightforward application of Lemma~\ref{LEMhardy-ineq-changevariables} yields sufficient conditions on the weights $u$ and $v$ for the inequality $\Vert   Tf\Vert_{L_q^u}\lesssim \Vert f\Vert_{L_p^v}$ to hold, where $T$ is an integral transform with splitting kernel. 

\begin{theorem}\label{THMsuff} 
	Let $T$ be an integral transform with splitting kernel, whose kernel satisfies \eqref{kernelest}. Let $1\leq p\leq q \leq \infty$ and let the weights $u,v$ be such that
	\begin{align}
		\sup_{r>0}\| w_1 u \|_{L_q(\varphi^{-1}(r,\infty))} \| s_1 v^{-1}\|_{L_{p'}(0,r)}&<\infty, \label{suffcond1}\\
		\sup_{r>0}\| w_2 u \|_{L_q(\varphi^{-1}(0,r))} \| s_2 v^{-1}\|_{L_{p'}(r,\infty)}&<\infty.\label{suffcond2}
	\end{align}
	Then the inequality $\Vert  Tf\Vert_{L_q^u}\lesssim \Vert f\Vert_{L_p^v}$ holds for every measurable $f$.
\end{theorem}
\begin{proof}
	First, it follows from \eqref{kernelest} that
	\begin{align*}
		\Vert  Tf\Vert_{L_q^u}&\lesssim \bigg(\int_0^\infty \big(u(y)w_1(y)I_{(0,\varphi(y))}(s_1 f)\big)^q dy\bigg)^{1/q}\nonumber \\
		&\phantom{=}+\bigg(\int_0^\infty \big(u(y)w_2(y)I_{(\varphi(y),\infty)}(s_2f)\big)^q dy\bigg)^{1/q}:=I_1+I_2.\label{EQpointwiseest}
	\end{align*}
	Now it suffices to apply Lemma~\ref{LEMhardy-ineq-changevariables} to deduce that
	\[
	I_1\lesssim \| f \|_{L_p^v}
	\]
	holds, provided that \eqref{suffcond1} is satisfied, and in a similar manner,
	\[
	I_2\lesssim  \| f \|_{L_p^v},
	\]
	provided that \eqref{suffcond2} holds.
\end{proof}

\subsection{Necessary conditions}
	
	We now obtain necessary conditions for \eqref{wni}. These are based in lower estimates for the kernel $K$ similar to \eqref{kernelest},  i.e., of the form
	\[
	K(x,y)\gtrsim \begin{cases}
		s_1(x)w_1(y),&x\in (0,\varphi(y)),\\
		s_2(x)w_2(y), &x\in (\varphi(y),\infty).
	\end{cases}
	\] 
	However, in contrast with sufficient conditions, the kernel $K$ may satisfy only one of these estimates, which will lead to partial necessary conditions. As it was already mentioned before, when $K$ satisfies both of these estimates, the necessary and sufficient conditions coincide (cf. Corollary~\ref{CORiff} below).

	We assume some technical conditions on the integrability (or boundedness) of the involved weights. More precisely, given a weighted Lebesgue space of the type $L_p^v$, with $1\leq p\leq \infty$, we will assume that
	\begin{equation}
		\label{EQintnearzero}
		s_1 v^{-1}\in L_{p'}(0,r), \qquad \text{for every }r>0,
	\end{equation}
	and
	\begin{equation}
		\label{EQintnearinfinity}
		s_2v^{-1}\in L_{p'}(r,\infty), \qquad  \text{for every }r>0.
	\end{equation}
	It is worth mentioning that \eqref{EQintnearzero} implies that $s_1v^{-1}\in L_1^{\textrm{loc}}$.
	
	The necessity result reads as follows.
	
	\begin{theorem}\label{THMnecessary} 
		Let $1\leq p,q\leq  \infty$, and consider the integral transform given by
        \[
            Tf(y) = \int_0^\infty f(x)K(x,y)\, dx.
        \]
        Let $\varphi:\R_+\to \R_+$ be a $C^1$ bijective function. Assume that the weighted norm inequality $\Vert  Tf\Vert_{L_q^u}\lesssim \Vert f\Vert_{L^v_p}$ holds for every $f\in L_p^v$.
		\begin{enumerate}
			\item If the kernel $K(x,y)$ satisfies the estimate
			\begin{equation}
				\label{EQlowerest1}
				K(x,y)\gtrsim s_1(x)w_1(y), \qquad x\in (0,\varphi(y)),
			\end{equation}
			with $s_1$ and $v$ satisfying \eqref{EQintnearzero}, then \eqref{suffcond1} holds.
			\item  If the kernel $K(x,y)$ satisfies the estimate
			\begin{equation}
				\label{EQlowerest2}
				K(x,y)\gtrsim s_2(x)w_2(y), \qquad x\in (\varphi(y),\infty),
			\end{equation}
			with $s_2$ and $v$ satisfying \eqref{EQintnearinfinity}, then \eqref{suffcond2} holds.
		\end{enumerate}
	\end{theorem}
	\begin{proof}
		We first prove the following claim:
		\begin{equation}
			\label{EQrelation}
			\text{if }r>0 \text{ and }x\in (0,r)\text{, then }x\in (0,\varphi(y)) \text{ for every }y\in \varphi^{-1}(r,\infty).
		\end{equation}
		We consider two different cases, according to whether $\varphi$ is decreasing or increasing. Firstly, if $\varphi$ is decreasing, then for $y\in \varphi^{-1}(r,\infty)=(0,\varphi^{-1}(r))$, we have
		\[
		y\leq \varphi^{-1}(r),\qquad \text{or equivalently}, \qquad r\leq \varphi(y),
		\]
		i.e., $(0,r)\subset (0,\varphi(y))$ for every $y\in \varphi^{-1}(r,\infty)$. Secondly, if $\varphi$ is increasing, then we have
		\[
		y\geq \varphi^{-1}(r),\qquad \text{or equivalently}, \qquad \varphi(y)\geq r,
		\]
            for all $y\in \varphi^{-1}(r,\infty)=(\varphi^{-1}(r),\infty)$.  This implies $(0,r)\subset (0,\varphi(y))$ for all $y\in \varphi^{-1}(r,\infty)$, which proves the claim. Similarly, one has that
		\begin{equation}
			\label{EQrelation2}
			\text{if }r>0 \text{ and }x\in (r,\infty)\text{, then }x\in (\varphi(y),\infty) \text{ for every }y\in \varphi^{-1}(0,r).
		\end{equation}
		We omit the details, which are analogous.
		
		We now subdivide the rest of the proof into three cases, namely, the case $1<p< \infty$, the case $p=\infty$, and the case $p=1$, which requires a slightly modified argument. 
		\begin{itemize}[leftmargin=*]
			\item Case $1<p<\infty$.
		\end{itemize}
		
		Let us prove the first part. We define, for $r>0$,
		\[
		f_r(x)=s_1(x)^{p'-1}v(x)^{-p'}\chi_{(0,r)}(x).
		\]
		On the one hand,
		\[
		\| f_r\|_{L_p^v} =\bigg(\int_0^r s_1(x)^{p'} v(x)^{p(1-p')}\,dx \bigg)^{\frac{1}{p}}=\| s_1 v^{-1}\|^{\frac{p'}{p}}_{L_{p'}(0,r)},
		\]
		which is finite, by \eqref{EQintnearzero}. On the other hand, by \eqref{EQlowerest1} and \eqref{EQrelation},
		\[
		Tf_r(y) =\int_0^r s_1(x)^{p'-1}v(x)^{-p'}K(x,y)\,dx\gtrsim w_1(y) \| s_1v^{-1}\|^{p'}_{L_{p'}(0,r)}, \qquad y\in \varphi^{-1}(r,\infty).
		\]
		Therefore,
		\begin{align*}
			\| s_1 v^{-1}\|^{\frac{p'}{p}}_{L_{p'}(0,r)}&=\| f_r\|_{L_p^v}\gtrsim \| Tf_r\|_{L_q^u} \geq \| Tf_r\|_{L_q^u(\varphi^{-1}(r,\infty))}\\
			&\gtrsim \| w_1u\|_{L_q(\varphi^{-1}(r,\infty))}\| s_1v^{-1}\|_{L_{p'}(0,r)}^{p'},
		\end{align*}
		i.e.,
		\[
		\| w_1u\|_{L_q(\varphi^{-1}(r,\infty))}\| s_1v^{-1}\|_{L_{p'}(0,r)}\lesssim 1.
		\]
		Taking the supremum on $r$ yields the desired result.

		The second part is analogous, and the details are omitted. In this case, it suffices to consider the function
		\[
		g_r(x)=s_2(x)^{p'-1} v(x)^{-p'}\chi_{(r,\infty)}(x),
		\]
		in place of $f_r$, and use a similar argument with \eqref{EQlowerest2} and \eqref{EQrelation2} in place of \eqref{EQlowerest1} and \eqref{EQrelation}, respectively.
		
		\begin{itemize}[leftmargin=*]
			\item Case $p=\infty$.
		\end{itemize}
		
		This case follows the same lines as the case $1<p<\infty$, with the difference that, in the first part, one should consider the function
		\[
		f_r(x)=v(x)^{-1}\chi_{(0,r)}(x),
		\]
		so that, similarly as above, since $Tf_r(y)\gtrsim w_1(y) \|s_1v^{-1}\|_{L_1(0,r)}$ for $ y\in \varphi^{-1}(r,\infty)$, we get
		\[
		1=\| f_r\|_{L_\infty^v}\gtrsim  \| Tf_r\|_{L_q^u(\varphi^{-1}(r,\infty))} \gtrsim \|w_1 u\|_{L_q(\varphi^{-1}(r,\infty))}\| s_1 v^{-1}\|_{L_{1}(0,r)},
		\]
		and the result follows by taking the supremum on $r$ (the last term is finite, by \eqref{EQintnearzero}). For the second part, we omit the details. It suffices to consider the function 
		\[
		g_r(x)=v(x)^{-1} \chi_{(r,\infty)}(x),
		\]
		and proceed similarly.
		\begin{itemize}[leftmargin=*]
			\item Case $p=1$.
		\end{itemize}
		
		Let, for $r>0$,
		\[
		f_r(x)=v(x)^{-1} h(x)\chi_{(0,r)}(x),
		\]
		where $h\in L^1(0,r)$ is nonnegative and different from the identically zero function. It is clear that $\| vf_r \|_{L_1(\R_+)}=\| h\|_{L_1(0,r)}$. Further,
		\[
		Tf_r(y)=\int_0^r v(x)^{-1} h(x) K(x,y)\,dx,
		\]
		for $y>0$. Hence,
		\[
		\| h\|_{L_1(0,r)}= \| vf_r \|_{L_1(\R_+)}\gtrsim \| uTf_r\|_{L_q(\R_+)} \geq \| uTf_r\|_{L_q(\varphi^{-1}(r,\infty))}.
		\]
		Now, by \eqref{EQlowerest1}, \eqref{EQrelation}, and the fact that $h$ is nonnnegative, we have 
		\[
		\| uTf_r\|_{L_q(\varphi^{-1}(r,\infty))} \gtrsim \|w_1 u\|_{L_q(\varphi^{-1}(r,\infty))}\| s_1 v^{-1}h\|_{L_{1}(0,r)},
		\]
		where the last term is finite, by \eqref{EQintnearzero} and the fact that $h\in L_1(0,r)$. Putting all estimates together, we obtain
		\[
		1\gtrsim \|w_1 u\|_{L_q(\varphi^{-1}(r,\infty))}\| s_1 v^{-1}h\|_{L_{1}(0,r)}\|h\|_{L_1(0,r)}^{-1},
		\]
		and replacing $h$ by $\widetilde{h}= h/\|h\|_{L_1(0,r)}$, we get
		\[
		1\gtrsim \|w_1 u\|_{L_q(\varphi^{-1}(r,\infty))} \sup_{\|\widetilde{h}\|_{L_1(0,r)}=1}\| s_1 v^{-1}\widetilde{h}\|_{L_{1}(0,r)}.
		\]
		Since $s_1$ and $v$ are nonnegative, the last supremum is equal to $\| s_1 v^{-1}\|_{L_\infty(0,r)}$ (see, e.g., \cite{BSbook}). Therefore, taking the supremum on $r$, we get
		\[
		1\gtrsim \sup_{r>0} \|w_1 u\|_{L_q(\varphi^{-1}(r,\infty))}\| s_1 v^{-1}\|_{L_\infty(0,r)},
		\]
		which completes the proof of the first part. 
		
		The proof of the second part is omitted, as it is analogous. In this case, one should consider
		\[
		g_r(x)=v(x)^{-1} h(x)\chi_{(r,\infty)}(x)
		\]
		in place of $f_r$ and proceed similarly.

	\end{proof}

	Theorems~\ref{THMsuff} and \ref{THMnecessary} allow us to characterize the weights $u$ and $v$ for which \eqref{wni} holds, provided that the involved kernel satisfies an asymptotic equivalence.
	\begin{corollary}\label{CORiff}
		Let $1<p\leq q\leq \infty$, and let the transform \eqref{trans} be such that
		\[
		K(x,y)\asymp \begin{cases}
			s_1(x)w_1(y),&\text{if }x\in (0,\varphi(y)),\\
			s_2(x)w_2(y), &\text{if }x\in (\varphi(y),\infty),
		\end{cases}
		\]
		where $\varphi$ is as in Theorem~\ref{THMsuff}, and $v$ and $s_j$, $j=1,2$, are such that \eqref{EQintnearzero} and \eqref{EQintnearinfinity} hold. Then, the weighted norm inequality $\Vert  Tf\Vert_{L_q^u}\lesssim \Vert  f\Vert_{L^v_p}$ is equivalent to the joint fulfillment of \eqref{suffcond1} and \eqref{suffcond2}.
	\end{corollary}
	
	As examples of integral transforms whose kernel satisfies the hypotheses of Corollary~\ref{CORiff} we can consider the Stieltjes transform, or the $\mathscr{H}_\alpha$ transform with $\alpha>-\frac{1}{2}$ (see Section~\ref{SECexamples} below).

	\subsection{A gluing lemma}
	It is interesting to note that under some assumptions on the functions $s_1, s_2, w_1$ and $w_2$, both conditions \eqref{suffcond1} and \eqref{suffcond2} (which, in essence, relate to Hardy-type inequalities) can be equivalently written as a single one. This is achieved through a so-called ``gluing lemma'' (cf. \cite{andersen,GKP}, where corresponding versions involving power weights were given in the case of classical Hardy inequalities). Here we present a full generalization, namely for Hardy inequalities involving changes of variables (i.e., Lemma~\ref{LEMhardy-ineq-changevariables}), and also with general weights.

	We have two versions of the gluing lemma, depending on whether the function $\psi=\varphi^{-1}$ in conditions \eqref{suffcond1} and \eqref{suffcond2} is increasing or decreasing. We start with the case where $\psi$ increasing.

	\begin{lemma}\label{lemmaglue}
		Let $f,g, s_j$, and $w_j$, $j=1,2$, be positive functions and $0<p,q<\infty$. Let $\psi:\R_+\to \R_+$ be an increasing bijective function. Assume that $s_2/s_1$ is nonincreasing, $w_1(\psi(t))\asymp s_2(t)$, and $w_2(\psi(t))\asymp s_1(t)$. Then, the conditions
		\begin{equation}\label{gluing1}
			\sup_{t>0} \bigg(\int_{\psi(t)}^\infty w_1(x)^qf(x)^q\, dx\bigg)^\frac{1}{q}  \bigg( \int_0^t s_1(x)^pg(x)^p\, dx\bigg)^\frac{1}{p}<\infty
		\end{equation}
		and
		\begin{equation}
			\label{gluing2}
			\sup_{t>0}\bigg( \int_{0}^{\psi(t)} w_2(x)^qf(x)^q\, dx\bigg)^\frac{1}{q}\bigg(\int_t^\infty s_2(x)^pg(x)^p\, dx\bigg)^\frac{1}{p}<\infty
		\end{equation}
		hold simultaneously if and only if 
		\begin{align}
			&\sup_{t>0} \bigg( w_1(\psi(t))^q \int_{0}^{\psi(t)} w_2(x)^qf(x)^q\, dx+ w_2(\psi(t))^q\int_{\psi(t)}^\infty w_1(x)^qf(x)^q\, dx\bigg)^\frac{1}{q}\nonumber \\
			&\times \bigg(\frac{1}{s_1(t)^p}\int_0^t s_1(x)^pg(x)^p\, dx+\frac{1}{s_2(t)^p}\int_t^\infty s_2(x)^pg(x)^p\, dx\bigg)^\frac{1}{p}<\infty.	\label{gluing3}
		\end{align}
	\end{lemma}
	\begin{proof}[Proof of Lemma~\ref{lemmaglue}]
		Since $w_1(\psi(t))\asymp s_2(t)$ and $w_2(\psi(t))\asymp s_1(t)$, it is clear that \eqref{gluing3} is equivalent to
		\begin{align}
			&\sup_{t>0}\bigg[ \frac{w_1(\psi(t))}{s_1(t)}\bigg(\int_{0}^{\psi(t)}w_2(x)^qf(x)^q\, dx\bigg)^\frac{1}{q}\bigg(\int_0^t s_1(x)^pg(x)^p\, dx\bigg)^\frac{1}{p} \nonumber\\
			&\phantom{=\sup}+\bigg(\int_{0}^{\psi(t)} w_2(x)^qf(x)^q\, dx\bigg)^\frac{1}{q}\bigg(\int_t^\infty s_2(x)^pg(x)^p\, dx\bigg)^\frac{1}{p} \nonumber\\
			&\phantom{=\sup}+\bigg(\int_{\psi(t)}^\infty w_1(x)^qf(x)^q\, dx\bigg)^\frac{1}{q} \bigg(\int_0^t s_1(x)^pg(x)^p\, dx\bigg)^\frac{1}{p} \nonumber\\
			&\phantom{=\sup}+\frac{w_2(\psi(t))}{s_2(t)}\bigg(\int_{\psi(t)}^\infty w_1(x)^qf(x)^q\, dx \bigg)^\frac{1}{q}\bigg(\int_t^\infty s_2(x)^pg(x)^p\, dx\bigg)^\frac{1}{p}\bigg]<\infty.	 \label{gluingsup}
		\end{align}		
		Thus, trivially, \eqref{gluing3} implies \eqref{gluing1} and \eqref{gluing2}. In order to prove the converse, we only have to show that the first and fourth terms of \eqref{gluingsup} are finite whenever \eqref{gluing1} and \eqref{gluing2} hold. We begin with the first term of \eqref{gluingsup}. For $t>0$, let $b(t)\in (0,t)$ be the number such that 
		\[
		\int_0^{b(t)} s_1(x)^pg(x)^p\, dx=\int_{b(t)}^t s_1(x)^pg(x)^p\, dx,
		\]
		which is well defined, given that all the involved functions are positive.
		Further, since $w_1(\psi(t))\asymp s_2(t)$, $w_2(\psi(t))\asymp s_1(t)$, $s_2/s_1$ is nonincreasing, and $\psi$ is increasing, then $w_1/w_2$ is almost decreasing, i.e., there exists a constant $C>0$ such that
		\[
		\frac{w_1(t)}{w_2(t)}\geq C\frac{w_1(t')}{w_2(t')},\qquad t'\geq t.
		\]
		Indeed, for $t'\geq t$, since $\psi$ is bijective and $\psi^{-1}$ is increasing,
		\[
		\frac{w_1(t)}{w_2(t)}\asymp \frac{s_2(\psi^{-1}(t))}{s_1(\psi^{-1}(t))}\geq \frac{s_2(\psi^{-1}(t'))}{s_1(\psi^{-1}(t'))}\asymp \frac{w_1(t')}{w_2(t')}.
		\]
		Thus,
		\begin{align*}
			&\phantom{=}\frac{w_1(\psi(t))}{s_1(t)}\bigg(\int_{0}^{\psi(t)}w_2(x)^qf(x)^q\, dx\bigg)^\frac{1}{q}\bigg(\int_0^t s_1(x)^pg(x)^p\, dx\bigg)^\frac{1}{p} \\ 
			&\asymp  \frac{w_1(\psi(t))}{s_1(t)} \bigg( \bigg(\int_0^{\psi(b(t))} w_2(x)^qf(x)^q\, dx\bigg)^\frac{1}{q} \bigg(\int_0^t s_1(x)^pg(x)^p\, dx\bigg)^\frac{1}{p}\\
			&\phantom{=}+\bigg(\int_{\psi(b(t))}^{\psi(t)}w_2(x)^qf(x)^q\, dx\bigg)^\frac{1}{q} \bigg(\int_{0}^t s_1(x)^pg(x)^p\, dx\bigg)^\frac{1}{p}\bigg)\\
			&\asymp \frac{s_2(t)}{s_1(t)} \bigg(\int_{0}^{\psi(b(t))}w_2(x)^qf(x)^q\, dx\bigg)^\frac{1}{q} \bigg(\int_{b(t)}^t s_1(x)^pg(x)^p\, dx\bigg)^\frac{1}{p}\\
			&\phantom{=}+\frac{w_1(\psi(t))}{w_2(\psi(t))}\bigg(\int_{\psi(b(t))}^{\psi(t)}w_2(x)^qf(x)^q\, dx\bigg)^\frac{1}{q} \bigg(\int_{0}^{b(t)} s_1(x)^pg(x)^p\, dx\bigg)^{\frac{1}{p}}\\
			&\lesssim  \bigg(\int_0^{\psi(b(t))} w_2(x)^qf(x)^q\, dx\bigg)^\frac{1}{q} \bigg(\int_{b(t)}^{t}s_2(x)^pg(x)^p\,dx\bigg)^\frac{1}{p}\\
			&\phantom{=}+\bigg(\int_{\psi(b(t))}^{\psi(t)}w_1(x)^qf(x)^q\, dx\bigg)^\frac{1}{q}\bigg(\int_0^{b(t)} s_1(x)^pg(x)^p\, dx\bigg)^\frac{1}{p}.
		\end{align*}
		The last expression is uniformly bounded in $t$, by \eqref{gluing1} and \eqref{gluing2}. We omit the estimate of the fourth term of \eqref{gluingsup}, since it essentially follows  the same steps as above. The only difference is that in place of $b(t)$, we choose the function $c(t)\in (t,\infty)$ so that
		\[
		\int_{\psi(t)}^{\psi(c(t))}w_2(x)^qf(x)^q\, dx=\int_{\psi(c(t))}^{\infty} w_2(x)^qf(x)^q\, dx, \qquad t>0.\qedhere
		\]
	\end{proof}
	
	A version of the gluing lemma for decreasing $\psi$ is stated without a proof, since it follows the same lines as that of Lemma~\ref{lemmaglue}. In this case, the corresponding Hardy-type conditions differ from \eqref{gluing1} and \eqref{gluing2}, in the sense that both integrals are taken in a neighbourhood of zero or of infinity, respectively (contrary to \eqref{gluing1} and \eqref{gluing2}, where, in each of the conditions, one integral is taken in a neighbourhood of zero, and the other one in a neighbourhood of infinity).
	\begin{lemma}\label{LEMMAgluingdual}
		Let $f,g, s_j$, and $w_j$, $j=1,2$, be positive functions and $0<p,q<\infty$. Let $\psi:\R_+\to \R_+$ be a decreasing bijective function. Assume that $s_2/s_1$ is nonincreasing, and $w_j(\psi(t))\asymp s_j(t)^{-1}$, $j=1,2$. Then, the conditions
		\begin{equation*}
			\sup_{t>0} \bigg(\int_{0}^{\psi(t)} w_1(x)^qf(x)^q\, dx\bigg)^\frac{1}{q}  \bigg( \int_0^t s_1(x)^pg(x)^p\, dx\bigg)^\frac{1}{p}<\infty
		\end{equation*}
		and
		\begin{equation*}
			\sup_{t>0}\bigg( \int_{\psi(t)}^{\infty} w_2(x)^qf(x)^q\, dx\bigg)^\frac{1}{q} \bigg(\int_t^\infty s_2(x)^pg(x)^p\, dx\bigg)^\frac{1}{p}<\infty
		\end{equation*}
		hold simultaneously if and only if 
		\begin{align*}
			&\sup_{t>0} \bigg(\frac{1}{w_1(\psi(t))^q} \int_{0}^{\psi(t)} w_1(x)^qf(x)^q\, dx+ \frac{1}{w_2(\psi(t))^q}\int_{\psi(t)}^\infty w_2(x)^qf(x)^q\, dx\bigg)^\frac{1}{q} \\
			&\times \bigg(\frac{1}{s_1(t)^p}\int_0^t s_1(x)^pg(x)^p\, dx+\frac{1}{s_2(t)^p}\int_t^\infty s_2(x)^pg(x)^p\, dx\bigg)^\frac{1}{p}<\infty.	
		\end{align*}
	\end{lemma}

	\medskip
	
	\section{Examples and comparison with known results}\label{SECexamples}
	\subsection{The sine transform}
	Here we restrict our attention to the case of power weights. As it was mentioned before, for the sine transform, the sufficient conditions from Theorem~\ref{THMsuff} are not sharp.
	
	Let us compare the sharp range for $\beta$ in the inequality $\Vert x^{-\beta}\widehat{f}_{\sin}\Vert_{L_q}\lesssim \Vert x^{\gamma}f\Vert_{L_p}$ with the range we obtain from Theorem~\ref{THMsuff}. It was recently proved in \cite{GLTPitt} that the inequality
	\begin{equation}
		\label{EQsinpitt}
		\Vert x^{-\beta}\widehat{f}_{\sin}\Vert_{L_q}\lesssim \Vert x^{\gamma}f\Vert_{L_p}
	\end{equation}
	holds if and only if
	\[
	\max\bigg\{0,\frac{1}{q}-\frac{1}{p'}\bigg\}\leq \beta<1+\frac{1}{q}, \qquad \beta=\gamma+\frac{1}{q}-\frac{1}{p'}.
	\]
	In our case, we can only deduce from Theorem~\ref{THMsuff} together with the estimate $|\sin xy|\leq \min\{xy,1\}$ that \eqref{EQsinpitt} holds if
	\[
	\frac{1}{q}< \beta<1+\frac{1}{q}, \qquad \beta=\gamma+\frac{1}{q}-\frac{1}{p'}.
	\]
	The fact that the sufficiency result does not yield a sharp range for $\beta$ comes from the fact that the estimate $|\sin xy|\leq 1$ is not optimal at all.
	
	In the necessity part, we can deduce from Theorem~\ref{THMnecessary} that \eqref{EQsinpitt} implies $\beta=\gamma+\frac{1}{q}-\frac{1}{p'}$ and $\beta<1+\frac{1}{q}$, which follows from the equivalence $\sin (xy) \asymp xy$ for $x\in (0,\frac{1}{y})$. In the interval $(\frac{1}{y},\infty)$ we cannot obtain any useful lower bound for $\sin( xy)$ that may be applicable in Theorem~\ref{THMnecessary}.
	
	\subsection{The Stieltjes transform}\label{SUBSECstieltjes}
	For the Stieltjes transform, Andersen proved in \cite{andersen} that the inequality $\Vert S_\lambda f\Vert_{L_q^u}\lesssim \Vert f\Vert_{L_p^v}$ holds if and only if
	\[
	\sup_{t>0}t^\lambda\bigg(\int_0^\infty \frac{u(x)^q}{(x+t)^{\lambda q}}\, dx\bigg)^{\frac{1}{q}}\bigg(\int_{0}^\infty \frac{v(x)^{-p'}}{(x+t)^{\lambda p'}}\, dx\bigg)^{\frac{1}{p'}}<\infty.
	\]
	We arrive to the same conclusion by applying Corollary~\ref{CORiff} together with the equivalence \eqref{EQstieltjesestimate} and the gluing lemma (Lemma~\ref{lemmaglue}). This result was also obtained in \cite[Proposition~4.6]{GKP} using the gluing lemma and the equivalence \eqref{EQstieltjesestimate}. 
	
	It is worth mentioning that Sinnamon characterized the inequality $\Vert S_\lambda f\Vert_{L_q^u}\lesssim \Vert f\Vert_{L_p^v}$ for the case $1\leq q<p\leq \infty$ in \cite{sinnamon}.
	
	\subsection{The \texorpdfstring{$\mathscr{H}_\alpha$} \text{} transform}
	For the $\mathscr{H}_\alpha$ transform and $1<p\leq q<\infty$, Rooney \cite{RonCan} found sufficient conditions on $\beta, \gamma$, so that the inequality
	\begin{equation}
		\label{EQHtranspitt}
		\Vert x^{-\beta} \mathscr{H}_\alpha f\Vert_{L_q}\lesssim \Vert x^\gamma f\Vert_{L_p}
	\end{equation}
	holds, namely  $\beta=\gamma+\frac{1}{q}-\frac{1}{p'}$ and
	\begin{align}
		\beta\geq \max\bigg\{ 0,\frac{1}{q}-\frac{1}{p'}\bigg\} \qquad \text{and}\qquad  \frac{1}{q}+\alpha-\frac{1}{2}<\beta&<\frac{1}{q}+\alpha+\frac{3}{2},& &\text{if }-\frac{1}{2}\leq \alpha<\frac{1}{2}, \nonumber \\
		\frac{1}{q}+\alpha-\frac{1}{2}<\beta&<\frac{1}{q}+\alpha+\frac{3}{2},& &\text{if }\alpha\geq \frac{1}{2}.\label{suffHtrans}
	\end{align}
	On the other hand, it was shown in \cite{miowni} that \eqref{suffHtrans} is also necessary for \eqref{EQHtranspitt} to hold in the case $\alpha>\frac{1}{2}$.
	
	Similarly as in the case of the sine transform, since the function $\mathbf{H}_\alpha$ is oscillating at infinity for $\alpha\leq \frac{1}{2}$ (cf. \eqref{eqEstimateStruveLarge}), we can only get sufficient conditions for \eqref{EQHtranspitt} that are not sharp by applying Theorem~\ref{THMsuff}. These read as $\beta=\gamma+\frac{1}{q}-\frac{1}{p'}$, and
	\[
	\frac{1}{q}<\beta<\frac{1}{q}+\alpha+\frac{3}{2}.
	\]
	In contrast, since for $\alpha>\frac{1}{2}$ we have the asymptotic equivalence 
	\begin{equation}\label{EQequivstruvefunct}
		(xy)^{\frac{1}{2}}\mathbf{H}_\alpha(xy)\asymp \begin{cases}
			(xy)^{\alpha+\frac{3}{2}},&\text{if }x\in (0,\frac{1}{y}),\\
			(xy)^{\alpha-\frac{1}{2}},&\text{if }x\in (\frac{1}{y},\infty),
		\end{cases}
	\end{equation}
	(cf. Remark~\ref{remarkHbigalpha}), we can characterize the weights $u$ and $v$ for which the inequality $\Vert \mathscr{H}_\alpha f\Vert_{L_q^u}\lesssim \Vert f\Vert_{L_p^v}$ holds.
	\begin{theorem}\label{THMiffstruve}
		Let $\alpha>\frac{1}{2}$ and $1\leq p\leq q\leq \infty$. Then, the inequality $\Vert \mathscr{H}_\alpha f\Vert_{L_q^u}\lesssim \Vert f\Vert_{L_p^v}$ holds if and only if the conditions 
		\begin{align}
			\sup_{t>0}\bigg( \int_{0}^{\frac{1}{t}} x^{(\alpha+\frac{3}{2})q}u(x)^q\, dx \bigg)^{\frac{1}{q}} \bigg(\int_0^t x^{(\alpha+\frac{3}{2})p'} v(x)^{-p'}\, dx\bigg)^{\frac{1}{p'}}&<\infty ,\label{EQsuffcondH1}\\
			\sup_{t>0}\bigg( \int_{\frac{1}{t}}^{\infty}x^{(\alpha-\frac{1}{2})q} u(x)^q\, dx\bigg)^{\frac{1}{q}}\bigg(\int_t^{\infty} x^{(\alpha-\frac{1}{2})p'} v(x)^{-p'}\, dx\bigg)^{\frac{1}{p'}}&<\infty,\label{EQsuffcondH2}
		\end{align}
		hold simultaneously.
	\end{theorem}
	\begin{proof}
		This follows by a direct application of Corollary~\ref{CORiff} and the equivalence \eqref{EQequivstruvefunct}.
	\end{proof}
	Theorem~\ref{THMiffstruve} improves the results from \cite{RonCan} in two directions, namely by considering general weights (thus also improving the results from \cite{miowni} in this respect) and by also giving necessary conditions for the corresponding inequality.
	
	We furthermore observe that we can use Lemma~\ref{LEMMAgluingdual} to rewrite both conditions \eqref{EQsuffcondH1} and \eqref{EQsuffcondH2} as a single one in the case $1<p\leq q<\infty$. Indeed, in this case $\psi(t)=\varphi^{-1}(t) =\frac{1}{t}$, $w_1(t)=s_1(t)=t^{\alpha+\frac{3}{2}}$, and $w_2(t)=s_2(t)=t^{\alpha-\frac{1}{2}}$. Thus, by Lemma~\ref{LEMMAgluingdual} (with $p'$ in place of $p$, $u(x)$ in place of $f(x)$, and $v(x)^{-1}$ in place of $g(x)$), conditions \eqref{EQsuffcondH1} and \eqref{EQsuffcondH2} are equivalent to the boundedness of 
	\begin{align*}
		\sup_{t>0}&\bigg(\frac{1}{t^{(-\alpha-\frac{3}{2})q}} \int_0^\frac{1}{t} x^{(\alpha+\frac{3}{2})q}u(x)^{q}\, dx + \frac{1}{t^{(-\alpha+\frac{1}{2})q}}  \int_\frac{1}{t}^\infty x^{(\alpha-\frac{1}{2})q} u(x)^{q}\, dx  \bigg)^\frac{1}{q}\\
		&\phantom{= \sup}\times \bigg(\frac{1}{t^{(\alpha+\frac{3}{2})p'}} \int_0^t x^{(\alpha+\frac{3}{2})p'}v(x)^{-p'}\, dx + \frac{1}{t^{(\alpha-\frac{1}{2})p'}}  \int_t^\infty x^{(\alpha-\frac{1}{2})p'} v(x)^{-p'}\, dx \bigg)^\frac{1}{p'}\\
		&=\sup_{t>0}t^{\alpha-\frac{1}{2}}\bigg(\int_0^\frac{1}{t} t^{2q} x^{(\alpha+\frac{3}{2})q}u(x)^{q}\, dx +  \int_\frac{1}{t}^\infty x^{(\alpha-\frac{1}{2})q} u(x)^{q}\, dx  \bigg)^\frac{1}{q}\\
		&\phantom{= \sup}\times  t^{-\alpha+\frac{1}{2}}\bigg( \int_0^t t^{-2p'}x^{(\alpha+\frac{3}{2})p'}v(x)^{-p'}\, dx +   \int_t^\infty x^{(\alpha-\frac{1}{2})p'} v(x)^{-p'}\, dx \bigg)^\frac{1}{p'}\\
		&=\sup_{t>0}\bigg(\int_0^\frac{1}{t} (xt)^{2q} x^{(\alpha-\frac{1}{2})q}u(x)^{q}\, dx +  \int_\frac{1}{t}^\infty x^{(\alpha-\frac{1}{2})q} u(x)^{q}\, dx  \bigg)^\frac{1}{q}\\
		&\phantom{= \sup}\times \bigg( \int_0^t \Big(\frac{x}{t}\Big)^{2p'}x^{(\alpha-\frac{1}{2})p'}v(x)^{-p'}\, dx +   \int_t^\infty x^{(\alpha-\frac{1}{2})p'} v(x)^{-p'}\, dx \bigg)^\frac{1}{p'}\\
		&=\sup_{t>0}\bigg(\int_0^\infty \min\big\{ (xt)^2,1\big\}^{q} x^{(\alpha-\frac{1}{2})q}u(x)^{q}\, dx \bigg)^\frac{1}{q} \\
            & \qquad\quad \times \bigg( \int_0^\infty \min\Big\{\Big(\frac{x}{t}\Big)^{2} ,1\Big\}^{p'} x^{(\alpha-\frac{1}{2})p'}v(x)^{-p'}\, dx \bigg)^\frac{1}{p'}.
	\end{align*}
	Taking into account that for $t,x>0$,
	\[
	\min\big\{ (xt)^2,1\big\} \asymp \frac{(xt)^2}{1+(xt)^2}=\frac{x^2}{t^{-2}+x^2}, \quad \textup{and} \quad \min\Big\{\Big(\frac{x}{t}\Big)^{2} ,1\Big\} \asymp \frac{\big(\frac{x}{t}\big)^{2}}{1+ \big(\frac{x}{t}\big)^{2}} = \frac{x^2}{t^2+x^2},
	\]
	we may rewrite Theorem~\ref{THMiffstruve} for the case $1<p\leq q<\infty$ as follows. 
	\begin{theorem}\label{CORgluningstruve}
		For $1<p\leq q<\infty$ and $\alpha>\frac{1}{2}$, the inequality $\Vert \mathscr{H}_\alpha f\Vert_{L_{q}^u}\lesssim \Vert f\Vert_{L_p^v}$ holds if and only if
		\[
		\sup_{t>0}  \bigg( \int_0^\infty \bigg(  \frac{x^{\alpha+\frac{3}{2}}}{t^{-2}+x^2}\bigg)^{q} u(x)^{q}\, dx\bigg)^\frac{1}{q}\bigg( \int_0^\infty \bigg(  \frac{x^{\alpha+\frac{3}{2}}}{t^{2}+x^2}\bigg)^{p'} v(x)^{-p'}\, dx\bigg)^\frac{1}{p'}<\infty.
		\]
	\end{theorem}

	\subsection{The Laplace transform}\label{SUBSEClaplace}
	\subsubsection{General weights}
	Weighted norm inequalities for the Laplace transform were first studied in \cite{BloomLaplace}. It was shown there that if
	\begin{equation}
		\label{EQsufflaplace}
		\sup_{t>0}\bigg(\int_0^{t} u(x)^q\, dx \bigg)^\frac{1}{q}\bigg(\int_0^{\frac{1}{t}} \bigg(\frac{1}{v}\bigg)^*(x)^{p'}\, dx\bigg)^\frac{1}{p'}<\infty,
	\end{equation}
	then the inequality 
	\begin{equation}
		\label{EQlaplace}
		\Vert \mathcal{L}f\Vert_{L_q^u}\lesssim \Vert f\Vert_{L_p^v},\qquad 1<p\leq q<\infty,
	\end{equation}
	holds (here $\big(\frac{1}{v}\big)^*$ denotes the decreasing rearrangement of $\frac{1}{v}$, cf. \cite{BSbook}). As for necessary conditions, it was proved that if \eqref{EQsufflaplace} holds, then
	\begin{equation}
		\label{EQneclaplace}
		\sup_{t>0}\bigg(\int_0^{t} u(x)^q\, dx \bigg)^{\frac{1}{q}}\bigg(\int_0^{\frac{1}{t}} v(x)^{-p'}\, dx\bigg)^{\frac{1}{p'}}<\infty.
	\end{equation}
	We note that in view of Theorem~\ref{THMnecessary}, such a condition is necessary for any integral transform whose kernel $K(x,y)$ satisfies an estimate of the form $K(x,y)\gtrsim 1$, for $x\in (0,\frac{1}{y})$. Whenever the weight $v$ is increasing, conditions \eqref{EQsufflaplace} and \eqref{EQneclaplace} coincide. On the other hand, if $v(x)=x^{\gamma}$ with $\gamma<0$, condition \eqref{EQsufflaplace} does not hold, since $\big(\frac{1}{v}\big)^*\equiv \infty$. For this reason, sufficient conditions not involving decreasing rearrangements were also given in \cite[Theorem~2]{BloomLaplace}. In particular, the condition
	\begin{equation}
		\label{EQnecsufstepanov}
		\sup_{t>0}\bigg(\int_0^\infty e^{-xtq} u(x)^q\, dx\bigg)^\frac{1}{q}\bigg(\int_0^tv(x)^{-p'}\, dx\bigg)^\frac{1}{p'}<\infty,
	\end{equation}
	is sufficient for \eqref{EQlaplace} to hold, whilst
	\begin{equation}
		\label{EQnecbloom}
		\sup_{t>0}\bigg(\int_0^\infty e^{-xt} u(x)^q\, dx\bigg)^\frac{1}{q}\bigg(\int_0^tv(x)^{-p'}\, dx\bigg)^\frac{1}{p'}<\infty,
	\end{equation}
	is necessary, and these conditions were shown to be equivalent whenever
    \[
        \int_{0}^{2y} u(x)^q\, dx\lesssim \int_0^y u(x)^q\, dx
        \qquad \textup{for} \quad y>0.
    \]
    Under this last assumption, condition \eqref{EQnecsufstepanov} (or \eqref{EQnecbloom}) was shown to characterize the boundedness of integral operators that generalize the Laplace transform in \cite[Theorem 3]{SU}. Corresponding results were obtained for the case $1<q<p<\infty$ both in \cite{BloomLaplace} and \cite{SU}, which we do not discuss here.

	Although conditions \eqref{EQsufflaplace} and \eqref{EQneclaplace} are close (and in fact equivalent if $v$ is nondecreasing), we observe that condition \eqref{EQsufflaplace} has a downside in the particular case of the Laplace transform. The kernel of the transform $K(x,y)=e^{-xy}$ decreases rapidly as $xy\to \infty$, which suggests that inequality \eqref{EQlaplace} could be satisfied with some power weight $v(x)=x^\gamma$, where $\gamma<0$, in which case $(1/v)^*(x)=\infty$ for $x>0$, and obviously \eqref{EQsufflaplace} does not hold. In fact, this is reflected in \eqref{EQnecsufstepanov}. We can also overcome this problem using Theorem~\ref{THMsuff}. 
	\begin{corollary}\label{CORlaplacesuff}
		Let $1\leq p\leq q\leq \infty$. The conditions
		\begin{align}
			\sup_{t>0}\| u \|_{L_q(0,t)} \| v^{-1}\|_{L_{p'}(0,\frac{1}{t})}&<\infty, \label{EQsuflap1}\\
			\sup_{t>0}\| x^{-n} u \|_{L_q(t,\infty)} \| x^{-n} v^{-1}\|_{L_{p'}(\frac{1}{t},\infty)}&<\infty,\quad \text{for some }n\in \N.\label{EQsuflap2}
		\end{align}
		are sufficient for the inequality $\| \mathcal{L}f\|_{L_q^u}\lesssim \| f\|_{L_p^v}$ to hold.
	\end{corollary}
	\begin{proof}
		It is enough use the estimate \eqref{EQestexponential} for the kernel of the Laplace transform and apply Theorem~\ref{THMsuff} with $s_1=w_1=1$, $s_2(t)=w_2(t)=t^{-n}$ (cf. \eqref{EQestexponential}), and $\varphi(y)=\frac{1}{y}$.
	\end{proof}
	
	\begin{remark}
		We note that for the case $1<p\leq q<\infty$, Corollary~\ref{CORlaplacesuff} is weaker than the sufficient condition \eqref{EQnecsufstepanov} given in \cite{BloomLaplace,SU}. More precisely, we can see that under the assumption  $1<p\leq q<\infty$, conditions \eqref{EQsuflap1} and \eqref{EQsuflap2} imply the sufficient condition \eqref{EQnecsufstepanov}. To see this, it suffices to use a similar argument as that of Corollary~\ref{CORgluningstruve} to write conditions \eqref{EQsuflap1} and \eqref{EQsuflap2} jointly. 
		
		Indeed, applying Lemma~\ref{LEMMAgluingdual} together with the estimate \eqref{EQestexponential} yields that the joint fulfillment of \eqref{EQsuflap1} and \eqref{EQsuflap2} is equivalent to 
		\begin{align*}
			&\phantom{=}\sup_{t>0}\bigg(\int_0^\infty \min\big\{1, (xt)^{-nq} \big\} u(x)^q\, dx\bigg)^\frac{1}{q}\bigg(\int_0^\infty \min\Big\{1,\Big(\frac{t}{x}\Big)^{np'} \Big\}v(x)^{-p'}\, dx \bigg)^{\frac{1}{p'}}\\
			&\asymp \sup_{t>0}\bigg(\int_0^\infty \Big( \frac{x^n}{t^{-n}+x^{n}}\Big)^q u(x)^q\, dx\bigg)^\frac{1}{q}\bigg(\int_0^\infty \Big( \frac{x^n}{t^{n}+x^{n}}\Big)^{p'}v(x)^{-p'}\, dx \bigg)^{\frac{1}{p'}}   <\infty,
		\end{align*}
		for some $n\in \N$. The estimate $e^{-xtq}\lesssim  \min\big\{1, (xt)^{-n} \big\}^q$ and the boundedness of the last supremum implies \eqref{EQnecsufstepanov}.
	\end{remark}
	

	\subsubsection{Power weights}
	Although Corollary~\ref{CORlaplacesuff} does not improve the sufficient condition given by \eqref{EQnecsufstepanov} in general, it still allows to characterize inequalities in the important case of power weights. 
	We first prove an auxiliary lemma.
	\begin{lemma}\label{LEMpowerweight}
		Let $u(x)=x^{-\beta}$ and $v(x)=x^\gamma$. For $1\leq p\leq q\leq \infty$, condition \eqref{EQsuflap1} implies \eqref{EQsuflap2}.
	\end{lemma}
	\begin{proof}
		We subdivide the proof in several cases.
		\begin{itemize}[leftmargin=*]
			\item Case $1<p\leq q<\infty$.
			We note that if \eqref{EQsuflap1} holds, then $\beta<\frac{1}{q}$ and $\gamma<\frac{1}{p'}$. Now,
			\[
			\| x^{-\beta}\|_{L_q(0,t)} \| x^{-\gamma}\|_{L_{p'}(0,\frac{1}{t})} \asymp t^{-\beta+\frac{1}{q}+\gamma-\frac{1}{p'}},
			\]
			and the latter is uniformly bounded if and only if $\beta =\gamma+\frac{1}{q}-\frac{1}{p'}$. On the other hand, for $n$ large enough, the norms in \eqref{EQsuflap2} are finite, and moreover
			\[
			\sup_{t>0}\| x^{-n-\beta}  \|_{L_q(t,\infty)} \| x^{-n-\gamma} \|_{L_{p'}(\frac{1}{t},\infty)} \asymp \sup_{t>0}t^{-\beta+\frac{1}{q}+\gamma-\frac{1}{p'}}<\infty.
			\]
			\item Case $1=p\leq q<\infty$. In this case, if \eqref{EQsuflap1} holds, then $\beta<\frac{1}{q}$ and $\gamma\leq 0$. Further,
			\[
			\| x^{-\beta}\|_{L_q(0,t)} \| x^{-\gamma}\|_{L_{\infty}(0,\frac{1}{t})} \asymp t^{-\beta+\frac{1}{q}+\gamma},
			\]
			which is uniformly bounded if and only if $\beta=\gamma+\frac{1}{q}$ (note that this condition combined with $\beta<\frac{1}{q}$ restricts $\gamma<0$ rather than $\gamma\leq 0$). As for \eqref{EQsuflap2}, if $n$ is large enough we get
			\[
			\sup_{t>0}\| x^{-n-\beta}  \|_{L_q(t,\infty)} \| x^{-n-\gamma} \|_{L_{\infty}(\frac{1}{t},\infty)} \asymp \sup_{t>0} t^{-\beta+\frac{1}{q}+\gamma}<\infty.
			\]
			\item Case $1<p\leq q=\infty$. The details are similar to the previous case and are omitted. We observe that if \eqref{EQsuflap1} holds, then necessarily $\beta\leq 0$, $\gamma<\frac{1}{p'}$, and $\beta=\gamma-\frac{1}{p'}$, which in turn implies that $\beta< 0$.
			\item Case $p=1$, $q=\infty$.  If \eqref{EQsuflap1} holds, then $\beta \leq 0$ and $\beta=\gamma$. Furthermore, \eqref{EQsuflap2} clearly holds. \qedhere
		\end{itemize}
	\end{proof}
	\begin{corollary}\label{CORlaplacepower}
		Let $1\leq p\leq q\leq \infty$, $(p,q)\neq (1,\infty)$. For $\beta,\gamma\in \R$, the inequality $\| x^{-\beta}\mathcal{L}f\|_{L_q}\lesssim \| x^\gamma f\|_{L_{p}}$ holds if and only if
		\[
		\beta<\frac{1}{q},\quad \text{and}\quad \beta =\gamma+\frac{1}{q}-\frac{1}{p'}, \quad \text{or equivalently,}\quad \gamma<\frac{1}{p'}\quad  \text{and}\quad \gamma=\beta-\frac{1}{q}+\frac{1}{p'}.
		\]
		The inequality $\| x^{-\beta}\mathcal{L}f\|_{L_\infty}\lesssim \| x^\gamma f\|_{L_{1}}$ holds if and only if $\beta=\gamma $ and $\beta\leq0$.
	\end{corollary}
	\begin{proof}
		By Corollary~\ref{CORlaplacesuff} and Lemma~\ref{LEMpowerweight}, condition \eqref{EQsuflap1} is sufficient for the desired inequality to hold. On the other hand, by Theorem~\ref{THMnecessary} and the estimate $K(x,y)=e^{-xy}\asymp 1$ for $x<\frac{1}{y}$, \eqref{EQsuflap1} is also necessary. Rewriting this condition in terms of the powers $\beta$ and $\gamma$ as in the proof of Lemma~\ref{LEMpowerweight} yields the conclusion.
	\end{proof}
	Note that in the case $1<p\leq q<\infty$, Corollary~\ref{CORlaplacepower} was obtained in \cite[Theorem~3.1]{BloomLaplace} (it also follows from \cite[Theorem 3]{SU}). We also stress that in the case of power weights, Corollary~\ref{CORlaplacepower} substantially improves the sufficient condition given by \eqref{EQsufflaplace}, which only allows to obtain a similar conclusion for the range $0\leq \gamma<\frac{1}{p'}$.

 \subsection*{Acknowledgments}
This research was supported by the Catalan AGAUR, through the Beatriu de Pin\'os program (2021BP  00072). The author is also grateful to the referee for useful comments that helped improving the quality of the paper, and especially Lemma~\ref{LEMpowerweight} and Corollary~\ref{CORlaplacepower}.


\EditInfo{December 30, 2023}{March 2, 2024}{Ana Cristina Freitas. Diogo Oliveira e Silva, Ivan Kaygorodov, and Carlos Florentino}

\end{document}